\documentstyle{article}
 
\newtheorem{theorem}{Theorem}[section]

\newenvironment{prooof}{\begin{description}
                   \item[{\small {\bf Proof:}}] \small}{\hfill {\bf Q.E.D.}
                                                          \medskip
                                                       \end{description}}
 
\newtheorem{defi}{Definition}[section]
\newtheorem{prop}{Proposition}[section]
\newtheorem{lemma}{Lemma}[section]
\newtheorem{rem}{Remark}[section]
\newtheorem{rems}{Remarks}[section]
\newtheorem{cor}{Corollary}[section]

\newcommand{\bdef}{\begin{defi}}
\newcommand{\ede}{\end{defi}}
\newcommand{\bsat}{\begin{theorem}}
\newcommand{\esat}{\end{theorem}}
\newcommand{\bprop}{\begin{prop}}
\newcommand{\eprop}{\end{prop}}
\newcommand{\blem}{\begin{lemma}}
\newcommand{\elem}{\end{lemma}}
\newcommand{\brem}{\begin{rem}}
\newcommand{\erem}{\end{rem}}
\newcommand{\brems}{\begin{rems}}
\newcommand{\erems}{\end{rems}}
\newcommand{\bcor}{\begin{cor}}
\newcommand{\ecor}{\end{cor}}
\newcommand{\bbew}{\begin{prooof}}
\newcommand{\ebew}{\end{prooof}}
\newcommand{\be}{\begin{equation}}
\newcommand{\ee}{\end{equation}}
\newcommand{\bea}{\begin{eqnarray}}
\newcommand{\eea}{\end{eqnarray}}
\newcommand{\beas}{\begin{eqnarray*}}
\newcommand{\eeas}{\end{eqnarray*}}
\newcommand{\ben}{\begin{enumerate}}
\newcommand{\een}{\end{enumerate}}
\newcommand{\lb}{\label}

\newcommand{\ra}{\rightarrow}

\newcommand{\f}{\frac}
\newcommand{\p}{\partial}

\newcommand{\Integer}{\:\mbox{\sf Z} \hspace{-0.82em} \mbox{\sf Z}\,}

\newcommand{\Real}{\mbox{I \hspace{-0.82em} R}}
\newcommand{\Complex}{
        \mbox{C \hspace{-1.16em} \raisebox{-0.018em}{\sf l}}\;}

\begin{document}
\title{On the deformation quantization of super-Poisson brackets 
       \vspace{1cm}}
\author{{\bf M.~Bordemann\thanks{mbor@phyq1.physik.uni-freiburg.de}} \\[3mm]
             Fakult\"at f\"ur Physik\\Universit\"at Freiburg\\
          Hermann-Herder-Str. 3\\79104 Freiburg i.~Br., F.~R.~G\\[3mm]
}
\date{FR-THEP-96/8 \\[1mm] May 1996 \\[2mm]
           }

%\begin{document}

\maketitle
\begin{abstract}
    We show that for every vector bundle $E$ over any given symplectic 
    manifold 
    there exists an explicitly given super Poisson bracket on the space of 
    sections of the
    dual Grassmann bundle associated to $E$ built out of the symplectic
    structure of $M$, a fibre metric on $E$ and a connection
    in $E$ compatible with the given fibre metric. Moreover, we construct
    a deformation quantization for this space of sections by means of
    a Fedosov type procedure.
\end{abstract}
\vfill
\newpage

\section {Introduction}

In the usual programme of deformation quantization (cf. \cite{BFFLS78})
the quantum mechanical multiplication is considered as a formal associative
deformation (a so-called star product)
of the pointwise multiplication of the classical observables, viz. the
algebra of smooth complex-valued functions on a given symplectic manifold.
The deformation is such that to first order in $\hbar$ the commutator of
the deformed product is proportional to the Poisson bracket. The difficult
question of existence of these star products for every symplectic manifold
was settled independently by DeWilde and Lecomte \cite{DL83} and
Fedosov \cite{Fed85}, \cite{Fed94}.

In the theory of supermanifolds the algebra ${\cal C}_0$ of classical 
superobservables
can be considered as the space of sections of the complexified 
dual Grassmann bundle
of an $n$-dimensional vector bundle $E$ (see e.~g. \cite{Bat80})
over a symplectic manifold $(M,\omega)$, i.e.
\be \lb{C0}
   {\cal C}_0:=\Complex\Gamma(\Lambda E^*), 
\ee
where the multiplication is the pointwise
wedge product. Clearly, ${\cal C}_0$ is a $\Integer_2$-graded
supercommutative algebra, i.e. $\phi\wedge\psi = (-1)^{d_1d_2}\psi\wedge\phi$
for $\phi,\psi\in\Gamma(\Lambda E^*)$, $\phi$ of degree $d_1$ and $\psi$ of
degree $d_2$. A {\em super-Poisson bracket} for ${\cal C}_0$ is by 
definition a $\Integer_2$-graded
bilinear map $M_1:{\cal C}_0\times {\cal C}_0\ra {\cal C}_0$ which is 
superanticommutative,
i.e. $M_1(\psi,\phi)=-(-1)^{d_1d_2}M_1(\phi,\psi)$, satisfies the 
superderivation
rule $M_1(\phi,\psi\wedge\chi)=M_1(\phi,\psi)\wedge\chi + 
(-1)^{d_1d_2}\psi\wedge M_1(\phi,\chi)$, and the super Jacobi identity, i.e.
$(-1)^{d_1d_3}M_1(M_1(\phi,\psi),\chi)+cycl.=0$ where $\chi\in {\cal C}_0$ is 
of degree $d_3$. \\
It is general not difficult 
to find super-Poisson brackets of purely algebraic type, i.e. which
vanish when one of their arguments is restricted to a smooth complex-valued
function, by means of a fibre metric $q$ in $E$ (see e.~g. \cite{BFFLS78},
p. 123, eqn 5-1):
\be
    M_1'(\phi,\psi)=q^{AB}(j(e_A)\phi)\wedge(i(e_B)\psi)
\ee
where $q^{AB}$ are the components of the induced fibre metric in the dual
bundle $E^*$ in the dual base to a local base $(e_A)$, $1\leq A\leq\dim E$,
of sections of
$E$, and $i(e_B)$ and $j(e_A)$ denote the usual interior product 
left antiderivation and right antiderivation, respectively. The definition
does not depend on the choice of that local base.\\
In case $M$ is $\Real^{2m}$ with the standard Poisson bracket 
one can combine the standard bracket with the above super-Poisson bracket 
to get
\be \lb{combined}
   M_1(\phi,\psi)= 
     \f{\p\phi}{\p q^i}\wedge\f{\p\psi}{\p p_i}
         -\f{\p\phi}{\p p_i}\wedge\f{\p\psi}{\p q^i}
          +~q^{AB}(j(e_A)\phi)\wedge(i(e_B)\psi) 
\ee
However, for nontrivial bundles it does not seem to be so 
obvious to generalize this bracket in the sense that it is
equal to -at least in degree zero- the Poisson bracket of the
base space $M$ when restricted to the sections of degree zero.

The main results of this paper are the following: firstly we construct an 
explicitly given
super-Poisson bracket $M_1$ for any vector bundle over any symplectic 
manifold which is equipped with an arbitrary nondegenerate fibre metric and a 
compatible connection
such that for two smooth complex-valued functions $f,g$ we have 
$M_1(f,g)=\{f,g\}~~+$ terms of higher Grassmann degree where
$\{~,~\}$ is the Poisson bracket on $M$. This superbracket generalizes
the above flat space superbracket (\ref{combined}), contains in a simple
polynomial manner the curvature of the connection in $E$, and does not yet
seem to have occurred in the literature as far as I know.\\
Secondly, we show by Fedosov's quantization procedure that the Grassmann
multiplication in ${\cal C}_0$ can formally be deformed into an associative
$\Integer_2$-graded multiplication $\ast$ on the space of formal power series 
in $\hbar$ with coefficients in ${\cal C}_0$,
\be \lb{C}
 {\cal C}:={\cal C}_0[[\hbar]],
\ee
such that the term proportional to $\hbar$ in that multiplication 
is equal to $({\bf i}/2)M_1$.

The paper is organized as follows:
We first transfer Fedosov's Weyl algebra bundle to our situation by
simply tensoring with the dual Grassmann bundle $\Lambda E^*$. The fibrewise
multiplication has also a component in $\Lambda E^*$ built by means of the
fibre metric in $E$. Then Fedosov's procedure can completely be imitated 
without further difficulties: we show the existence of a Fedosov connection
$D$ of square zero whose kernel in the space of antisymmetric degree zero
is in linear $1-1$ correspondence to $\cal C$ which immediately gives rise
to the desired quantum deformation (Theorem \ref{ass}).\\
In the third part we explicitly compute the super-Poisson bracket $M_1$
as the term proportional to $({\bf i}\hbar)/2$ by means of Fedosov's 
recursion formulae (Theorem \ref{thebracket}).

Notation: In all of this paper the Einstein sum convention is used that two
equal indices are automatically summed up over their range unless stated
otherwise. Moreover, we widely make use of Fedosov's notation in 
\cite{Fed94} with the following exceptions: we use the symbol $\nabla$
to denote the covariant derivative and not Fedosov's $\p$ and 
describe the occurring
symmetric tensor fields with $\vee$ products (see e.g. \cite{Gre78}, 
p. 209-226) and use the symmetric substitution operator $i_s$ instead of
Fedosov's functions of $y$ and derivatives with respect to $y$. 

\section{The Fedosov quantization procedure for dual Grassmann bundles}

Let $(M,\omega)$ be a $2m$-dimensional symplectic manifold and $E$ an 
$n$-dimensional real vector bundle
over $M$ with a fixed positive definite 
fibre metric $q$. For the computations that will follow we shall use
co-ordinates $(x^1,\cdots,x^{2m})$ in a chart $U$ of $M$. The base fields
$\f{\p}{\p x^i}$ will be denoted by $\p_i$ ($1\leq i\leq 2m)$ for short.
For computations in $E$ we shall use a local base $(e_A)$, $(1\leq A\leq n)$
of sections of $E$ over $U$. Denote the dual base in the dual bundle $E^*$ 
of $E$ by $(e^A)$, $(1\leq A\leq n)$. Let $\Lambda\in\Gamma(\Lambda^2TM)$
denote the Poisson structure of $(M,\omega)$, i.e. the Poisson bracket of two
smooth real valued functions $f,g$ is given by $\{f,g\}:=\Lambda(df,dg)$.
Denoting the components of $\omega$ and $\Lambda$ in that chart by
$\omega_{ij}:=\omega(\p_i,\p_j)$ and $\Lambda^{ij}:=\Lambda(dx^i,dx^j)$
we use the sign conventions of \cite{AM85} where 
$\Lambda^{ik}\omega_{jk}=\delta^i_j$. Fix a torsion-free symplectic
connection $\nabla^M$ in the tangent bundle of $M$. This is well-known to
always exist which can be seen by He\ss 's formula 
$\omega(\nabla^M_XY,Z):=\omega(\tilde{\nabla}_XY,Z)
+\f{1}{3}(\tilde{\nabla}_X\omega)(Y,Z)+\f{1}{3}(\tilde{\nabla}_Y\omega)(X,Z)$
where $X,Y,Z$ are arbitrary vector fields on $M$ and $\tilde{\nabla}$ is
an arbitrary torsion-free connection on $M$ (see \cite{Hes81}). 
Fix a connection $\nabla^E$
in $E$ which is compatible with $q$, i.e. $X(q(e_1,e_2))=q(\nabla^E_Xe_1,e_2)
+q(e_1,\nabla^E_Xe_2)$ for an arbitrary vector field $X$ on $M$ and sections
$e_1,e_2$ of $E$.

We are now forming the {\em Fedosov algebra} ${\cal W}\otimes\Lambda$:
\be \lb{FedAlg}
  {\cal W}\otimes\Lambda:=\big( \times_{s=0}^\infty
  \Gamma(\Complex(\vee^sT^*M\otimes\Lambda E^*\otimes \Lambda T^*M))\big)
      [[\hbar]]
\ee
This is to say that the elements of ${\cal W}\otimes\Lambda$ are formal
sums $\sum_{s,t=0}^\infty w_{st}\hbar^t$ where the $w_{st}$ are smooth
sections in the complexification of the vector bundle 
$\vee^sT^*M\otimes\Lambda E^*\otimes \Lambda T^*M$. In what follows we shall 
sometimes use the following factorized sections 
$F:=f\otimes \phi\otimes \alpha\hbar^{t_1}$ and 
$G:=g\otimes \psi\otimes\beta\hbar^{t_2}$ where $f\in\Gamma(\vee^{s_1}T^*M)$, 
$g\in \Gamma(\vee^{s_2}T^*M)$, $\phi\in\Gamma(\Lambda^{d_1}E^*)$, 
$\psi\in\Gamma(\Lambda^{d_2}E^*)$, $\alpha\in\Gamma(\Lambda^{a_1}T^*M)$, and
$\beta\in\Gamma(\Lambda^{a_2}T^*M)$. Let $deg_s,deg_E,deg_a,deg_{\hbar}$ be
the obvious degree maps from ${\cal W}\otimes\Lambda$ to itself, i.e.
those $\Complex$-linear maps for which the above factorized sections
$f\otimes \phi\otimes \alpha\hbar^{t_1}$ are eigenvectors to the eigenvalues
$s_1,d_1,a_1,t_1$ respectively and which we refer to as the symmetric degree,
the $E$-degree, the antisymmetric degree, and the $\hbar$-degree, 
respectively. Moreover, let $P_E$ and $P_{\hbar}$ be the corresponding
maps $(-1)^{deg_E}$ and $(-1)^{deg_{\hbar}}$ which we refer to as the
$E$-parity and the $\hbar$-parity, respectively. We say that a 
$\Complex$-linear endomorphism $\Phi$ of ${\cal W}\otimes\Lambda$
is of $\zeta$-degree $k$ ($\zeta=s,a,E,\hbar$) iff 
$[deg_{\zeta},\Phi]=k\Phi$. Analogously, $\Phi$ is said to be of 
$\zeta$-parity $(-1)^k$ ($\zeta=E,\hbar$) iff 
$P_{\zeta}\Phi P_{\zeta}=(-1)^k\Phi$.
Let $C$ denote the
complex conjugation of sections in ${\cal W}\otimes\Lambda$.\\
We shall sometimes write $\cal W$ for the space of elements of 
${\cal W}\otimes\Lambda$ having zero antisymmetric degree and
${\cal W}\otimes\Lambda^a$ for the space of those elements having 
antisymmetric degree $a$. The space ${\cal W}\otimes\Lambda$ is an 
associative algebra with respect to the usual pointwise product where
we do {\em not} use the graded tensor product of the two Grassmann algebras
involved. More precisely, for the above factorized sections the pointwise
or undeformed multiplication is simply given by
\be \lb{undefmulti}
  (f\otimes \phi\otimes \alpha\hbar^{t_1})
    (g\otimes \psi\otimes\beta\hbar^{t_2})
     :=(f\vee g)\otimes(\phi\wedge\psi)\otimes(\alpha\wedge\beta)
             \hbar^{t_1+t_2}.
\ee
Note that the above four degree maps are derivations 
and the above two parity maps are automorphisms of this multiplication.
Moreover, ${\cal W}\otimes\Lambda$ is supercommutative in the sense that
\be \lb{supercomm}
     GF=(-1)^{d_1d_2+a_1a_2}FG
\ee
A linear endomorphism $\Phi$ of ${\cal W}\otimes\Lambda$
of $E$-parity $(-1)^{d'}$ and antisymmetric degree $a'$ is said to be a
superderivation of type $((-1)^{d'},a')$ of the undeformed algebra
${\cal W}\otimes\Lambda$ iff $\Phi(FG)=(\Phi F)G+(-1)^{d'd_1+a'a_1}F(\Phi G)$.
Let $\sigma$ denote the linear map
\be \lb{sigma}
  \sigma:{\cal W}\otimes\Lambda\ra 
             \Gamma(\Lambda E^*\otimes \Lambda T^*M)[[\hbar]]
\ee
which projects onto the component of symmetric degree zero and clearly is
a homomorphism for the undeformed multiplication.\\
We now combine the two covariant derivatives $\nabla^M_X$ in $TM$ and
$\nabla^E_X$ in $E$ into a covariant derivative $\nabla_X$ in $TM\otimes E$
in the usual fashion and extend it canonically to a connection $\nabla$
in ${\cal W}\otimes\Lambda$. On the above factorized sections we get in
a chart:
\be \lb{connection}
   \nabla(f\otimes\phi\otimes\alpha)
      = \big( (\nabla^M_{\p_i}f)\otimes\phi
          + f\otimes(\nabla^E_{\p_i}\phi)\big) \otimes(dx^i\wedge\alpha)
            + f\otimes\phi\otimes d\alpha.
\ee
\\
In order to define a deformed fibrewise associative multiplication consider
the following insertion maps for a vector field $X$ on $M$ and a section
$e$ of $E$: let $i_a(X)$ and $i(e)$ denote the usual inner product 
antiderivations in $\Gamma(\Lambda T^*M)$ and $\Gamma(\Lambda E^*)$, 
respectively, and extend them in a canonical manner to superderivations of
type $(1,-1)$ and $(-1,0)$ of the undeformed algebra 
${\cal W}\otimes\Lambda$, respectively.
Let $j(e)$ be defined by $P_Ei(e)$. Moreover,
let $i_s(X)$ denote the corresponding inner product derivation (or symmetric 
substitution, \cite{Gre78}, p.209-226) in  
$\times_{s=0}^\infty\Gamma(\vee^sT^*M)$, again extended to a derivation of
the undeformed algebra
${\cal W}\otimes\Lambda$ in the canonical way. Let $q^{AB}$ denote the 
components
of the induced fibre metric $q^{-1}$ in $E^*$, i.e. 
$q^{AB}:=q^{-1}(e^A,e^B)$. Note that $q^{AB}$ is the inverse matrix
to $q(e_A,e_B)$. Then for two elements $F,G$ of ${\cal W}\otimes\Lambda$ we 
can now define the fibrewise deformed multiplication $\circ$:
\bea
    F\circ G & := & \sum_{k,l=0}^\infty \f{(i\hbar/2)^{k+l}}{k!l!}
                                             \nonumber \\
             &    & \big( \Lambda^{i_1j_1}\cdots\Lambda^{i_kj_k}
                     (i_s(\p_{i_1})\cdots i_s(\p_{i_k})F)
                      (i_s(\p_{j_1})\cdots i_s(\p_{j_k})G) \nonumber \\
             &   &      +~q^{A_1B_1}\cdots q^{A_lB_l}
                             (j(e_{A_1})\cdots j(e_{A_l})F)
                               (i(e_{B_1})\cdots i(e_{B_l})G)  \big) 
                                                 \lb{superFed}
\eea
Moreover, let $\delta$ and $\delta^*$ be the canonical superderivations
of the undeformed algebra ${\cal W}\otimes\Lambda$ of type $(1,1)$ and 
$(1,-1)$, respectively, which are induced by the identity map of $T^*M$ to
$T^*M$ where in the case of $\delta$ the preimage of the identity is regarded
as being part of $\vee T^*M$ and the image as being part of $\Lambda T^*M$,
and vice versa in the case of $\delta^*$. On the above factorized sections
these maps read in co-ordinates
\bea
    \delta(f\otimes\phi\otimes\alpha) &  = &
                       (i_s(\p_i)f)\otimes\phi\otimes(dx^i\wedge\alpha) 
                                      \lb{delta} \\
    \delta^*(f\otimes\phi\otimes\alpha) &  = &
                    (dx^i\vee f)\otimes\phi\otimes(i_a(\p_i)\alpha). 
                                      \lb{deltastern}
\eea
Define the total degree $Deg$:
\be \lb{Deg}
   Deg:=2deg_{\hbar}+deg_s+deg_E
\ee
A $\circ$-superderivation of type $((-1)^{d'},a')$ is defined in an analogous
manner as for the undeformed multiplication.\\
We collect some properties of the above structures in the following
\bprop \lb{eprop} 
  With the above definitions and notations we have the following:
 \ben
  \item $\delta^2=0=(\delta^*)^2$ and 
         $\delta\delta^*+\delta^*\delta=deg_s+deg_a$.
  \item $\delta\nabla+\nabla\delta=0$.
  \item $Ker(\delta)\cap Ker(deg_a)={\cal C}$.
  \item $P_E$ is a
    $\circ$-auto\-mor\-phism and $deg_a$ is a $\circ$-derivation which
    equips the Fedosov algebra $({\cal W}\otimes\Lambda,\circ)$ with the 
    structure of a $\Integer_2\times\Integer$-graded associative algebra.
  \item $\delta$, $\nabla$, and $Deg$ are $\circ$-superderivations
    of type $(1,1)$, $(1,1)$, and $(1,0)$, respectively.
  \item The parity map $P_{\hbar}$ and the complex conjugation $C$
    are graded $\circ$-anti\-au\-to\-mor\-phisms, i.e.
     $\Phi(F\circ G)=(-1)^{d_1d_2+a_1a_2}G\circ F$ for $\Phi=P_{\hbar}, C$.
 \een
\eprop
\bbew
1. Straight forward.\\
2. This follows from the vanishing torsion of 
$\nabla^M$.\\ 
3. Without the factor $\Lambda E^*$ the kernel of $\delta$ in the space of
  antisymmetric degree zero consists of the constants, which proves this 
  statement.\\
4. The associativity of $\circ$ is known, see
e.g. \cite{BFFLS78}, p. 123, eqn 5-2, and can be done by a long straight 
forward computation. We shall sketch a shorter proof: $\circ$ is defined
on each fibre (for $m\in M$) 
${\cal W}_m:=(\times_{i=0}^\infty (\vee^i T_mM^*\otimes \Lambda E_m^*\otimes
\Lambda T_mM^*))[[\hbar]]$ on which we can rewrite the multiplication in the 
more compact form ($F,G\in{\cal W}_m$)
\bea
   F\circ G & = & \mu(
     e^{\f{{\bf i}\hbar}{2}(\Lambda^{ij}_mi_s(\p_i)\otimes i_s(\p_j)
         +q^{AB}_mj(e_A)\otimes i(e_B))}(F\otimes G)) \nonumber \\
            & =: & 
           \mu(e^{\f{{\bf i}\hbar}{2}(R+S)}(F\otimes G))
                                         \lb{compactcirc}
\eea
where the tensor product is over $\Complex[[\hbar]]$ and 
$\mu$ denotes the undeformed fibrewise multiplication. Due to the derivation
properties of $i_s(\p_i)$, $i(e_A)$, and  $j(e_B)$ we get formulas like
\beas 
     R~\mu\otimes 1  & = & \mu\otimes 1~(R_{13}+R_{23})  \\
     R~1\otimes \mu  & = & 1\otimes \mu~(R_{12}+R_{13})  \\
     S~\mu\otimes 1  & = & \mu\otimes 1~(S_{13}(P_E)_2+S_{23}) \\
     S~1\otimes\mu   & = & 1\otimes\mu~(S_{12}+S_{13}(P_E)_2)
\eeas
where the index notation is borrowed from Hopf algebras and indicates
on which of the three tensor factors of ${\cal W}_m$ the maps $R$, $S$,
and $P_E$ should act, e.g. $R_{23}:=1\otimes R$, 
$(P_E)_2:=1\otimes P_E\otimes 1$. These ``pull through formulas'' can
be used to pull through the corresponding formal exponentials. Since
all the maps $i_s(\p_i)$ commute with $j(e_A)$ and $i(e_B)$ and since the
$j(e_A)$ commute with all $i(e_B)$ whereas $j(e_A)$ and $j(e_B)$ anticommute
as well as $i(e_A)$ and $i(e_B)$ we can conclude that all the six 
maps $R_{12}$, $R_{13}$, $R_{23}$, $S_{12}$, $S_{13}(P_E)_2$, and $S_{23}$ 
pairwise commute. This is the essential step
for associativity. The gradation properties are immediate.\\
5. The derivation properties of $\delta$ and $Deg$ are clear, for the 
corresponding statement for $\nabla$ the fact that $\nabla^M$ preserves
the Poisson structure $\Lambda$ and that $\nabla^E$ preserves the dual
fibre metric $q^{-1}$ is crucial.\\
6. Straight forward.
\ebew

Due to the first part of this proposition we can construct a 
$\Complex[[\hbar]]$-linear endomorphism $\delta^{-1}$ of the Fedosov
algebra in the following way:
on the above factorized sections $F$ we put
\be \lb{deltaminus1}
 \delta^{-1}F := \left\{
 \begin{array}{cr}
     \f{1}{s_1+a_1}\delta^*F & {\rm ~if~}s_1+a_1\geq 1
                                                      \\
         0                   & {\rm ~if~}s_1+a_1=0
 \end{array} 
  \right.
\ee

Since ${\cal W}\otimes\Lambda$ is an $\Integer_2\times\Integer$-graded
associative algebra we can form the $\Integer_2\times\Integer$-graded
super Lie bracket which reads on the above factorized sections:
\be \lb{supercom}
   [F,G]:=ad(F)G:=F\circ G-(-1)^{d_1d_2+a_1a_2}G\circ F
\ee
It follows from the associativity of $\circ$ that $ad(F)$ is 
$\circ$-superderivation of the Fedosov algebra 
$({\cal W}\otimes\Lambda,\circ)$ of type $((-1)^{d_1},a_1)$. Note that
the map $\f{{\bf i}}{\hbar}ad(F)$ which we shall often use in what follows
is always well-defined because of the
supercommutativity of the undeformed multiplication (\ref{supercomm}).\\
Consider now the curvature tensors $R^M$ of $\nabla^M$ and $R^E$ of 
$\nabla^E$, i.e. for three vector fields $X,Y,Z$ on $M$ and a section $e$ 
of $E$ we have 
$R^M(X,Y)Z=\nabla^M_X\nabla^M_YZ-\nabla^M_Y\nabla^M_XZ-\nabla^M_{[X,Y]}Z$
and 
$R^E(X,Y)e=\nabla^E_X\nabla^E_Ye-\nabla^E_Y\nabla^E_Xe-\nabla^E_{[X,Y]}e$.
Define elements $R^{(M)}$ and $R^{(E)}$ of the Fedosov algebra which are 
contained in
$\Gamma(\vee^2T^*M\otimes\Lambda^2T^*M)$ and 
$\Gamma(\Lambda^2E^*\otimes\Lambda^2T^*M)$, respectively, as follows
where $V,W$ are vector fields on $M$ and $e_1,e_2$ are sections of $E$:
\bea
   R^{(M)}(V,W,X,Y)     & := & \omega(V,R^M(X,Y)W) \lb{RM} \\
   R^{(E)}(e_1,e_2,X,Y) & := & -q(e_1,R^E(X,Y)e_2). \lb{RE}
\eea
Note that this is well-defined: since $\nabla^M$ preserves $\omega$ and
$\nabla^E$ preserves $q$ it follows that $R^{(M)}$ is symmetric in $V,W$
and $R^{(E)}$ is antisymmetric in $e_1,e_2$. In co-ordinates these
two elements of the Fedosov algebra can be written in the form
$R^{(M)}=(1/4)R^{(M)}_{klij}dx^k\vee dx^l\otimes 1\otimes dx^i\wedge dx^j$
and 
$R^{(E)}=(1/4)R^{(E)}_{ABij}1\otimes e^A\wedge e^B\otimes dx^i\wedge dx^j$.
Set
\be \lb{Fedcurv}
   R:=R^{(M)}+R^{(E)}.
\ee
Then the following Proposition is immediate:
\bprop \lb{rprop}
 With the above definitions and notations we have:
 \ben
  \item $\nabla^2=\f{{\bf i}}{\hbar}ad(R)$.
  \item $P_E(R)=R$, $P_{\hbar}(R)=R$ and $C(R)=R$.
  \item $\delta R=0$.
  \item $\nabla R=0$.
 \een
\eprop
\bbew
  1. Straightforward computation. \\
  2. Obvious.\\
  3. This is a consequence of the vanishing torsion of $\nabla^M$
     (first Bianchi identity). \\
  4. This is a reformulation of the second Bianchi identity for linear
     connections in arbitrary vector bundles.
\ebew
We shall now make the ansatz for a Fedosov connection $D$, i.e. we are 
looking for an element $r\in{\cal W}\otimes\Lambda^1$ {\em of even E-parity},
i.e. $P_E(r)=r$, such that the map
\be \lb{Fedder}
    D:=-\delta+\nabla+\f{{\bf i}}{\hbar}ad(r)
\ee
has square zero, i.e. $D^2=0$. The following properties of $D$ for any $r$
are crucial:
\blem \lb{D2} 
 Let $r$ be an arbitrary element of ${\cal W}\otimes\Lambda^1$ of
 even $E$-parity. Then
  \ben
   \item $D^2=\f{{\bf i}}{\hbar}ad(-\delta r+\nabla r+R
                   +\f{{\bf i}}{\hbar}r\circ r)$.
   \item $D(-\delta r+\nabla r+R+\f{{\bf i}}{\hbar}r\circ r)=0$.
  \een 
\elem
\bbew 
 This is straight forward using Proposition \ref{rprop} and 
 the fact that $r\circ r=\f{1}{2}[r,r]$ for the above elements $r$ of 
 even $E$-parity and odd antisymmetric degree.
\ebew

For an arbitrary element $w\in{\cal W}\otimes\Lambda$ we shall make the
following decomposition according to the total degree $Deg$:
\be
     w=\sum_{k=0}^\infty w^{(k)}~~~{\rm where~}Deg(w^{(k)})=kw^{(k)}
\ee
Note that each $w^{(k)}$ is always a {\em finite} sum of sections in
some $\Gamma(\vee^sT^*M\otimes \Lambda E^*\otimes \Lambda T^*M)$. The 
subspaces of all elements of ${\cal W}\otimes\Lambda$, ${\cal W}$,
${\cal W}\otimes\Lambda^a$, and ${\cal C}$
of total degree $k$ will be denoted
by ${\cal W}^{(k)}\otimes\Lambda$, ${\cal W}^{(k)}$,
${\cal W}^{(k)}\otimes\Lambda^a$, and ${\cal C}^{(k)}$,
respectively.

As in Fedosov's paper \cite{Fed94} there is the following 
\bsat \lb{rexist}
 With the above definitions and notations: Let  
 $r\in{\cal W}\otimes\Lambda^1$ be defined by the following recursion:
 \beas
       r^{(3)}   & := & \delta^{-1}R \\
       r^{(k+3)} & := & \delta^{-1}\left( \nabla r^{(k+2)}+\f{{\bf i}}{\hbar}
                          \sum_{l=1}^{k-1}r^{(l+2)}\circ r^{(k-l+2)}
                                            \right)
 \eeas
 Then $r$ has the following properties: it is real ($C(r)=r$),
 depends only on $\hbar^2$ ($P_{\hbar}(r)=r)$, has even $E$-parity, and
 is in the kernel of $\delta^{-1}$.\\
 Moreover, the corresponding Fedosov
 derivation $D=-\delta+\nabla+({\bf i}/\hbar)ad(r)$ has square zero.
\esat
\bbew The behaviour of $r$ under the parity transformations and complex
 conjugation immediately follows from the fact that they commute with
 $\delta^{-1}$ and from their (anti)homomorphism properties
 (Prop.\ref{eprop}, 3., 5.; Prop.\ref{rprop}, 2.).\\
 Let $A:=-\delta r+\nabla r+R+\f{{\bf i}}{\hbar}r\circ r=:-\delta r+R+B$.
 Recall the equation $\delta\delta^{-1}+\delta^{-1}\delta=1$ on the
 subspace of the Fedosov algebra where $deg_s+deg_a$ have nonzero 
 eigenvalues. Clearly,
 $A^{(2)}=-\delta r^{(3)}+R=0$ because $\delta R=0$ (Prop.\ref{rprop},3.)
 hence $R=\delta\delta^{-1}R$.
 Suppose $A^{(l)}=0$ for all $2\leq l\leq k+1$. By Lemma \ref{D2}, 2.
 we have $0=(DA)^{(k+1)}=
 -\delta A^{(k+2)}=-\delta B^{(k+2)}$. Hence 
 $B^{(k+2)}=\delta\delta^{-1}B^{(k+2)}=\delta r^{(k+3)}$ proving
 $A^{k+2}=0$ which inductively implies $D^2=0$ since we had already shown
 that $r$ is of even $E$-parity.
\ebew

We shall now compute the kernel of the Fedosov derivation. More precisely,
define
\be \lb{WD}
  {\cal W}_D:=Ker(D)\cap Ker(deg_a).
\ee
As in Fedosov's paper \cite{Fed94} we have the important characterization:
\bsat \lb{theoWD}
 With the above definitions and notations: ${\cal W}_D$ is a subalgebra
 of the Fedosov algebra $({\cal W}\otimes\Lambda,\circ)$. Moreover, the
 map $\sigma$ (\ref{sigma}) restricted to ${\cal W}_D$ 
 is a $\Complex[[\hbar]]$-linear bijection onto 
 ${\cal C}$.
\esat
\bbew
 The kernel of a superderivation is always a subalgebra. Since $D$ and 
 $\sigma$ are $\Complex[[\hbar]]$-linear the subalgebra ${\cal W}_D$
 is a $\Complex[[\hbar]]$-submodule of ${\cal W}$.\\
 Let 
 $w\in{\cal W}$. Decompose $w=w_0+w_+$ where $w_0:=\sigma(w)$ and
 $w_+:=(1-\sigma)(w)$. We shall prove by induction over the total degree $k$
 that $w\in{\cal W}$ is in ${\cal W}_D$ iff for all nonnegative integers
 $k$ $w^{(k)}_0$ is arbitrary
 in ${\cal C}^{(k)}$ and $w^{(k)}_1$ is uniquely given by the equation
 \be \lb{w+}
    w^{(k)}_+=\delta^{-1}\left( \nabla w^{(k-1)}+
                \f{{\bf i}}{\hbar}\sum_{l=1}^{k-2}\left[ 
               r^{(l+2)}, w^{(k-1-l)}\right]
                             \right)=:(Aw)^{(k)}
 \ee
 where of course an empty sum is defined to be zero and $w^{(0)}_1=0$.
 Note that $Dw=-\delta w+Aw$ and that the $\Complex$-linear map $A$ does not
 lower the total degree of $w$.\\
 Now the equation $(Dw)^{(k)}=0$ is equivalent to the inhomogeneous
 equation $\delta w^{(k+1)}=(Aw)^{(k)}$. A necessary condition for this
 equation to be solvable for $w^{(k+1)}$ clearly is 
 $\delta((Aw)^{(k)})=0$. But this is also
 sufficient since then $(Aw)^{(k)}=\delta\delta^{-1}(Aw)^{(k)}$ and we
 have the particular solution $w^{(k+1)}_+=\delta^{-1}(Aw)^{(k)}$
 (since $\sigma\delta^{-1}=0$) which satisfies (\ref{w+}). To this
 particular solution any solution to the homogeneous equation 
 $\delta w^{(k)}=0$ can be added which precisely is the space 
 ${\cal C}^{(k)}$.\\
 It remains to show that conversely every initial piece
 $w':=w^{(0)}_0+w^{(1)}_0+w^{(1)}_++\cdots+w^{(k)}_0+w^{(k)}_+$
 where $w^{(l)}_0$ was arbitrarily chosen in ${\cal C}^{(l)}$,
 $w^{(l)}_+$ is determined by (\ref{w+}) for all $0\leq l\leq k$, and
 $(Dw')^{(l)}=0$ for all $-1\leq l\leq k-1$ can be continued to
 $w'':=w'+w^{(k+1)}_0+w^{(k+1)}_+$ with $w^{(k+1)}_0$ arbitrary in
 ${\cal C}^{(k+1)}$, $w^{(k+1)}_+$ determined by (\ref{w+}), and
 $(Dw'')^{(k)}=0$. By induction, this will eventually lead to 
 $w\in{\cal W}_D$
 characterized by the above properties. Indeed, since $D^2=0$
 we have $0=(D^2w')^{(k-1)}=-\delta((Dw')^{(k)})=-\delta((Aw')^{(k)})$.
 Define $w^{(k+1)}_+$ by $\delta^{-1}((Aw')^{(k)})$ 
 and choose any $w^{(k+1)}_0\in{\cal C}^{(k+1)}$.
 It follows at once that $w^{(k+1)}_+$ satisfies (\ref{w+}) and that
 we get $(Dw'')^{(k)}=0$ which proves the induction and the Theorem. 
\ebew
Let
\be \lb{tau}
  \tau : {\cal C}\ra {\cal W}_D\subset {\cal W}
\ee
be the inverse of the restriction of $\sigma$ to ${\cal W}_D$. For
$\phi\in\Gamma(\Lambda E^*)$ we shall speak of $\tau(\phi)$ as the
{\em Fedosov-Taylor series of} $\phi$ and refer to the components
$\tau(\phi)^{(k)}$ as the Fedosov-Taylor coefficients. 
We collect some of the properties of $\tau$ in the following
\bprop \lb{tautheo}
  With the above definitions and notations:
 \ben
  \item $\tau$ commutes with $P_E$, $P_{\hbar}$, and $C$.
  \item Let $\phi=\sum_{d=0}^n\phi^{(d)}\in\Gamma(\Lambda E^*)$ where
    $n:=\dim E$. 
    Then
    $Deg(\phi^{(d)})=d\phi^{(d)}=deg_E(\phi^{(d)})$.\\
    Moreover
    \bea
        \tau(\phi)^{(0)} & = & \phi^{(0)} \\
        \tau(\phi)^{(1)} & = & \delta^{-1}(\nabla\phi^{(0)})+\phi^{(1)} \\
         \vdots          &   &        \vdots \nonumber \\
       \tau(\phi)^{(n)}  & = & \delta^{-1}\left(
                                    \nabla(\tau(\phi)^{(n-1)})
                                 +\f{{\bf i}}{\hbar}\sum_{l=1}^{n-2}
                           \left[ r^{(l+2)},\tau(\phi)^{(n-1-l)}
                               \right] \right) 
                                      +\phi^{(n)} \nonumber \\
                         &   &               \\
      \tau(\phi)^{(n+1)} & = & \delta^{-1}\left(
                                    \nabla(\tau(\phi)^{(n)})
                                 +\f{{\bf i}}{\hbar}\sum_{l=1}^{n-1}
                           \left[ r^{(l+2)},\tau(\phi)^{(n-l)}\right]
                                                 \right) \\
          \vdots         &   &   \vdots  \nonumber \\
      \tau(\phi)^{(k+1)} & = & \delta^{-1}\left(
                                    \nabla(\tau(\phi)^{(k)})
                                 +\f{{\bf i}}{\hbar}\sum_{l=1}^{k-1}
                         \left[ r^{(l+2)},\tau(\phi)^{(n-l)}\right]
                                                 \right)
    \eea
    where $k\geq n$.
    The Fedosov-Taylor series $\tau(\phi)$ depends only on $\hbar^2$.
    \item For any nonnegative integer $k$ the map 
    $\phi\mapsto\tau(\phi)^{(k)}$ is a polynomial in $\hbar$
    whose coefficients are differential operators from 
    $\Gamma(\Lambda E^*)$ into some $\Gamma(\vee^sT^*M\otimes\Lambda E^*)$
    of order $k$.
 \een
\eprop
\bbew
 Since $r$ is invariant under the parity maps and complex conjugation,
 it follows that $D$ commutes with these three maps, hence ${\cal W}_D$
 is stable under these maps. Since $\sigma$ obviously commute with them,
 so does the inverse of its restriction to ${\cal W}_D$, $\tau$. 
 The rest is a consequence of the preceding Theorem and a straight
 forward induction.
\ebew

Define the following $\Complex[[\hbar]]$-bilinear multiplication on
${\cal C}$: for $\phi,\psi\in{\cal C}$ 
\be \lb{star}
     \phi\ast\psi:=\sigma(\tau(\phi)\circ\tau(\psi)).
\ee
We shall call $\ast$ the {\em Fedosov star product associated to}
$(M,\omega,\nabla^M,E,q,\nabla^E)$.
For $\phi,\psi\in\Gamma(\Lambda E^*)$ the star product $\phi\ast\psi$
will be a formal power series in $\hbar$ which we shall write in the 
following form:
\be \lb{dieMs}
   \phi\ast\psi =: \sum_{t=0}^\infty \left( \f{{\bf i}\hbar}{2} \right)^t
                                 M_t(\phi,\psi).   
\ee
We list some important properties of the Fedosov star product in the
following
\bsat \lb{ass} With the above definitions and notations:
 \ben
  \item The Fedosov star product is associative and $\Integer_2$-graded,
    i.e. $P_E$ is an automorphism of $({\cal C},\ast)$. The map 
    $P_{\hbar}$ and the complex
    conjugation $C$ are graded antiautomorphisms of $({\cal C},\ast)$.
  \item The $\Complex$-bilinear maps $M_t$ are all bidifferential, real,
    vanish on the constant functions in each argument for $t\geq 1$, and
    have the following symmetry property:
    \be \lb{symprop}
         M_t(\psi,\phi)=(-1)^{t}(-1)^{d_1d_2}M_t(\phi,\psi).
    \ee
  \item The term of order $0$ is equal to the pointwise Grassmann 
    multiplication. Hence $({\cal C},\ast)$ is a formal associative
    deformation of the supercommutative algebra 
    $({\cal C}_0,\wedge)$.
 \een
\esat
\bbew Basically, every stated property is easily derived from the 
 definitions (\ref{star}) and (\ref{dieMs}) and the corresponding behaviour
 of the fibrewise multiplication $\circ$ under $P_E$, $P_\hbar$, and $C$.
 The reality of the $M_t$ follows easily from the graded antihomomorphism
 property of $C$ once eqn (\ref{symprop}) is proved by means of the graded
 antihomomorphism property of the $\hbar$-parity. Since $\tau(1)$ is easily
 seen to be equal to $1$ we have $1\ast\psi=\psi=\psi\ast 1$, and the $M_t$
 must vanish on $1$ for $t\geq 1$. Finally, each $M_t$ obviously depends
 on only a finite number of Fedosov-Taylor coefficients whence it must
 be bidifferential.
\ebew

\section{Computation of the super-Poisson bracket}

In this section we are going to compute an explicit expression for the
term $M_1$ of the Fedosov star product defined in the last section
(compare (\ref{dieMs}) and Theorem \ref{ass}). Only by means of the
graded associativity of the deformed algebra $({\cal C},\ast)$ we can
derive the following
\blem \lb{SuperPoisson}
 Let $\phi,\psi,\chi$ be sections in ${\cal C}_0$ of $E$-degree
 $d_1,d_2,d_3$, respectively. Then
 \bea
  M_1(\psi,\phi)  & = & -(-1)^{d_1d_2}M_1(\phi,\psi) \\
  M_1(\phi,\psi\wedge\chi) & = & M_1(\phi,\psi)\wedge\chi
                                   +(-1)^{d_1d_2}\psi\wedge 
                                          M_1(\phi,\chi) \\
  0 & = & (-1)^{d_1d_3}M_1(M_1(\phi,\psi),\chi)+{\rm cycl.} 
 \eea
 Hence $M_1$ is a {\em super-Poisson bracket} on ${\cal C}_0$.
\elem
\bbew The first property is a particular case of (\ref{symprop}).
 Consider now the graded commutator 
 $[\phi,\psi]:=\phi\ast\psi-(-1)^{d_1d_2}\psi\ast\phi$ on $\cal C$.
 Because of the graded associativity of $\ast$ we have the superderivation
 property $[\phi,\psi\ast\chi]=[\phi,\psi]\ast\chi+
 (-1)^{d_1d_2}\psi\ast [\phi,\chi]$. Writing this out with the $M_t$
 and taking the term of order $\hbar$ we get the second property.
 For the third, take the term of order $\hbar^2$ in the super Jacobi identity
 for the graded commutator.
\ebew
Before we are going to compute $M_1$ directly it is useful to introduce 
the following notions:\\
For $\phi$ in ${\cal C}_0$ let $\phi_1$ and $\rho$ denote the 
component of symmetric degree one and $\hbar$-degree zero of the 
Fedosov-Taylor coefficient $\tau(\phi)$ and the section $r$ 
(Theorem \ref{rexist}), respectively. Note that $\phi_1$ is a smooth
section
in the bundle $T^*M\otimes \Lambda E^*$. Denote by $\Lambda_0 E^*$ the 
subbundle of the dual Grassmann bundle consisting of elements of even
degree. Then $\rho$ is a smooth
section in $T^*M\otimes \Lambda_0 E^*\otimes T^*M$. Consider now
the bundle $TM\otimes \Lambda_0 E^*\otimes T^*M$. There is an obvious
fibrewise associative multiplication $\bullet$ in that bundle which comes 
from the identification of $TM\otimes T^*M$ with the bundle of linear 
endomorphism of $TM$: let $X,Y$ be vector fields on $M$, 
$\phi,\psi\in\Lambda_0 E^*$, and $\alpha,\beta$ one-forms on $M$. Then
\be \lb{bullet}
     (X\otimes\phi\otimes\alpha)\bullet(Y\otimes\psi\otimes\beta)
        :=(\alpha(Y))X\otimes(\phi\wedge\psi)\otimes\beta.
\ee
Let $\hat{R}^E$
be the section in $\Gamma(TM\otimes\Lambda^2 E^*\otimes T^*M)$ whose
components in a bundle chart read
\be \lb{hatre}
  \hat{R}^E:=\f{1}{4}\Lambda^{ik}R^{(E)}_{ABkj}
                \p_i\otimes e^A\wedge e^B \otimes dx^j
                  =:\p_i\otimes (\hat{R}^E)^i_j\otimes dx^j,
\ee
and let $\hat{\rho}\in\Gamma(TM\otimes\Lambda_0 E^*\otimes T^*M)$ be defined
by
\be \lb{hatrho}
  \hat{\rho}:=\p_i \otimes\Lambda^{ik}i_s(\p_k)\rho
             =:\p_i\otimes \hat{\rho}^i_j\otimes dx^j.
\ee
Note that we can 
form arbitrary power series in $\hat{R}^E$ by using the multiplication 
$\bullet$ since $\hat{R}^E$ is nilpotent.
 
We have the following 
\blem With the above notations and definitions:
 \bea
  M_1(\phi,\psi) & = & \Lambda^{ij}(i_s(\p_i)\phi_1)(i_s(\p_j)\psi_1) 
                           +q^{AB}(j(e_A)(\phi))(i(e_B)(\psi)) \\
      \phi_1     & = & dx^j((1-\hat{\rho})^{-1})^i_j\nabla^E_{\p_i}\phi \\
  \hat{\rho}     & = & 1-(1-2\hat{R}^E)^{1/2}.
 \eea
 where $(1-\hat{\rho})^{-1}$ and $(1-2\hat{R}^E)^{1/2}$ denote the 
 corresponding power series with respect to the $\bullet$ multiplication.
\elem
\bbew
 The first equation is a straight forward computation.\\
 For the second, use the Fedosov recursion for $\tau(\phi)$, 
 (Proposition \ref{tautheo}),
 note that $\phi_1^{(k)}$ is zero for $k\geq n+2$ and that only the
 component $\rho$ of $r$ matters since both $\tau(\phi)$ and $r$ depend
 only on $\hbar^2$, sum over the total degree which yields the equation
 \[ 
  \phi_1 = \delta^{-1}\nabla^E\phi+dx^j(\hat{\rho})^i_j(i_s(\p_i)\phi_1) 
 \]
 which proves the second equation.\\
 For the third, use the Fedosov recursion for $r$, (Theorem \ref{rexist}), 
 take the
 component of symmetric degree 1 and $\hbar$-degree zero, sum over the total
 degree, and arrive at the quadratic equation
 \[
    \hat{\rho}-\hat{R}^E=\f{1}{2}\hat{\rho}\bullet\hat{\rho}.
 \]
 Since $r$ and hence $\rho$ does not contain components of symmetric degree
 zero, there is only one solution to this equation, namely the above third 
 equation.
\ebew

This Lemma immediately implies the desired formula for the super-Poisson
bracket:
\bsat \lb{thebracket}
   The super-Poisson bracket $M_1$ obtained by the Fedosov star product
   takes the following form:
   \beas
     M_1(\phi,\psi) & = & 
           \Lambda^{ij}((1-2\hat{R}^E)^{-1/2})^k_i\wedge
                         ((1-2\hat{R}^E)^{-1/2})^l_j\wedge
                            \nabla^E_{\p_k}\phi\wedge\nabla^E_{\p_l}\psi \\
           & & +~q^{AB}(j(e_A)(\phi))(i(e_B)(\psi))
   \eeas
\esat
\bbew 
   Clear from the Lemma !
\ebew

\brems
 \ben
  \item In case $(M,\omega)$ is K\"ahler there exist star
   products of Wick type on $M$ (see \cite{BW96}): they are charactrized
   by the property 
   that for any two complex-valued smooth functions $f,g$
   on $M$ the star product $f\ast'g$ is made out of bidifferential 
   operators which differentiate $f$ in holomorphic
   directions only and $g$ in antiholomorphic directions only. It seems
   to me very likely that super analogues of these star products can
   readily be formulated for any complex holomorphic hermitean vector
   bundle over $M$.
  \item If the dual Grassmann bundle $\Lambda E^*$ is replaced by the
   symmetric power $\vee E^*$ and the fibre metric $q$ by some antisymmetric
   bilinear form on the fibres covariantly constant by some connection
   in $E$ the whole construction can presumably carried 
   through as well.
 \een
\erems

\vspace{0.5cm}

\noindent {\Large Acknowledgment}

\vspace{0.2cm}

\noindent I would like to thank R.~Eckel for useful discussions, 
 S.~Waldmann for a critical reading of the manuscript, E.~Galois
 for constant encouragement, and J.~S.~B. for his cool stuff in minor,
 in particular BWV 889 and BWV 893.

\end{document}